\documentclass[11pt]{amsart}
\usepackage{epsfig}
\usepackage{amssymb}
\usepackage{amscd}
\usepackage{epsf}

\newcommand{\sevafig}[3]{\begin{figure}[h]\centerline{
 \epsfig{file=#1,width=#2,angle=#3}}
\bigskip\caption{}\end{figure}}

\newcommand{\nw}{\bigwedge\nolimits}

\renewcommand{\phi}{\varphi}
\newcommand{\suml}{\sum\limits}

\newcommand{\ndot}{\bullet}

\def\matho#1{\mathop{\mathrm{#1}}}

\def\suml{\sum\limits}
\DeclareMathSizes{11.1}{10}{8}{6}

\newcommand{\Lie}{{\matho{Lie}}}

\newcommand{\Hom}{\matho{Hom}\nolimits}

\newcommand{\Aut}{\matho{Aut}}

\newcommand{\Hoch}{\mathrm{Hoch}}
\newcommand{\Clif}{\mathrm{Clif}}

\newcommand{\C}{\mathbb C}

\newtheorem*{theorem}{Theorem}

\newtheorem*{lemma}{Lemma}
\newtheorem*{corollary}{Corollary}
\newtheorem*{conjecture}{Conjecture}

\theoremstyle{remark}
\newtheorem*{remark}{Remark}

\theoremstyle{definition}

\author{Boris Shoikhet}
\title
{Gerstenhaber bracket on double Hochschild complex and deformation theory}
\date{07.10.1999}
\address{IUM, 11 Bol'shoj Vlas'evskij per.,
Moscow 121002, Russia}
\email{borya@mccme.ru}

\begin{document}

\maketitle

\begin{abstract}

We construct a differential and a Lie bracket on the space\linebreak
$\{\Hom (A^{\otimes k}, A^{\otimes l})\},_{k,l\ge 0}$ for any associative 
algebra
$A$. The restriction of this bracket to the space $\{\Hom (A^{\otimes k},
A)\},_{k\ge 0}$ is exactly the Gerstenhaber bracket. We discuss some
formality conjecture related with this construction. We also
discuss some applications to deformation theory.

\end{abstract}

\setcounter{subsection}{-1}

\subsection{}
There exists a well-known way (due to Jim Stasheff) to define the
Hochschild differential and the Gerstenhaber bracket via
coderivations on the cofree coalgebra cogenerated by the vector
space $A[1]$ ($A$ is an associative algebra). Analogously, one
can define a dual differential and bracket on the graded space
$\{\Hom (A, A^{\otimes k})\},_{k\ge 0}$ using derivations of tensor
algebra, generated by the space $A[1]$ for any coalgebra $A$. In
both cases the bracket does not depend on the (co)algebra
structure, and only the differential does.

In these notes we generalize this construction involving all
``differential operators'' on the tensor algebra $T^{\ndot}(A[1])$,
not only of the first order. In such a way, we define a
bidifferential and a bracket on the bigraded space $\{\Hom
(A^{\otimes k}, A^{\otimes l})\},_{k,l\ge 0}$. It turns out that
for any associative algebra $A$ \emph{with unit} the total
cohomology of the bicomplex are equal examples to zero. So, the
interesting examples appear for the algebras \emph{without unit},
for example, for the algebras of polynomials without unit.

\subsection{}

Let $V$ be a vector space, and let $\Psi_1\colon V^{\otimes k_1}\to
V^{\otimes l_1}$,\ \ $\Psi_2\colon V^{\otimes k_2}\to V^{\otimes l_2}$ be
any two maps. We are going to define the bracket $[\Psi_1,\Psi_2]$.

Let $\Psi\colon V^{\otimes k}\to V^{\otimes l}$ be any map. It defines a map
$i(\Psi)\colon V^{\otimes N}\to V^{\otimes(N-k+l)}$,\ \ $N\gg0$, as follows:
\begin{multline}
i(\Psi)(v_1\otimes v_2\otimes\ldots\otimes v_n)=\Psi(v_1\otimes\ldots\otimes
v_k)\otimes v_{k+1}\otimes\ldots\otimes v_N+\\
+v_1\otimes\Psi(v_2\otimes\ldots\otimes v_{k+1})\otimes\ldots\otimes
v_N+\ldots+\\
+v_1\otimes v_2\otimes\ldots\otimes\Psi(v_{N-k+1}\otimes\ldots\otimes v_N).
\end{multline}

It is clear that the composition $i(\Psi_1)\circ i(\Psi_2)$ has \emph{not}
the
form $i(\Psi_3)$ for some~$\Psi_3$. Nevertheless, one has the following
statement.

\begin{lemma}
The bracket $i(\Psi_1)\circ i(\Psi_2)-i(\Psi_2)\circ i(\Psi_1)$ has a form
$i\left(\suml_{s\ge0}\Psi_{3,s}\right)$ where $\Psi_{3,s}$ is a map
$V^{\otimes(k_1+k_2-s)}\to V^{\otimes(l_1+l_2-s)}$.
\end{lemma}

\begin{proof}
It is clear.
\end{proof}

We set:
\begin{equation}
[\Psi_1,\Psi_2]=\sum_s\Psi_{3,s}.
\end{equation}

\begin{remark}
The case $s=0$ may appear only when $k_1=k_2=0$ or $l_1=l_2=0$. In the
general case, $s\ge1$.
\end{remark}

We have constructed a ``Lie algebra of differential operators'' on the
tensor algebra $T^\ndot(V)$. Let us note that this Lie algebra is not
corresponded to an associative algebra.

\subsection{}

Let $V=A[1]$ where $A$ is an associative algebra, $m_A\colon A^{\otimes2}\to
A$ is the product. Let $\Psi\colon A^{\otimes k}\to A^{\otimes l}$ be any
map.

\begin{lemma}
The bracket $[m_A,\Psi]=\Psi_1+\Psi_2$, where $\Psi_1\colon
A^{\otimes(k+1)}\to A^{\otimes l}$ and \hbox{$\Psi_2\colon A^{\otimes k}\to
A^{\otimes(l-1)}$} are defined as follows\emph:
\begin{multline}
\Psi_1(a_1\otimes\ldots\otimes
a_{k+1})=(a_1\otimes1\otimes1\otimes\dots)\cdot\Psi(a_2\otimes\ldots\otimes
a_{k+1})-\\
-\Psi(a_1a_2\otimes a_3\otimes\ldots\otimes a_{k+1})+\Psi(a_1\otimes
a_2a_3\otimes\dots)\pm\\
\pm\Psi(a_1\otimes\ldots\otimes a_k)\cdot(1\otimes1\otimes\ldots\otimes
a_{k+1});
\end{multline}
\begin{equation}
\Psi_2(a_1\otimes\ldots\otimes a_k)=\pm(m_{12}-m_{23}+\ldots\pm m_{l-1,l})\circ
\Psi(a_1\otimes\ldots\otimes a_k),
\end{equation}
where
\begin{equation}
m_{i,i+1}(b_1\otimes\ldots\otimes b_l)=b_1\otimes\ldots\otimes
b_{i-1}\otimes b_i\cdot b_{i+1}\otimes b_{i+2}\otimes\ldots\otimes b_l.
\end{equation}
\end{lemma}

\begin{proof}
It is a direct calculation.
\end{proof}

\begin{corollary}
$$
[m_A,m_A]=0
$$
\end{corollary}

\qed

Let us denote $\Psi_1(\Psi)=d_1(\Psi)$, $\Psi_2(\Psi)=d_2(\Psi)$ (in the
notations of Lemma above); it follows from Corollary that $d_1+d_2$ is a
bidifferential, i.\,e. $d_1^2=0$,\  \ $d_2^2=0$,\ \ $d_1d_2=\pm d_2d_1$. We
have the following bicomplex:
\begin{gather*}
\begin{CD}
0 @>>>\C @>0>> A^* @>d_1>>(A^{\otimes2})^* @>d_1>>(A^{\otimes3})^*
@>d_1>>{}\\
@. @A0AA @AA0A @AA0A @AA0A \\
0 @>>> A @>d_1>> \Hom(A,A) @>d_1>> \Hom(A^{\otimes2},A) @>d_1>>
\Hom(A^{\otimes3},A) @>d_1>> \relax \\
@. @Ad_2AA @AAd_2A @AAd_2A @AAd_2A \\
0 @>>> A^{\otimes2} @>d_1>> \Hom(A,A^{\otimes2}) @>d_1>>
\Hom(A^{\otimes2},A^{\otimes2}) @>d_1>>
\Hom(A^{\otimes3},A^{\otimes2}) @>d_1>> \relax
\\
@. @Ad_2AA @AAd_2A @AAd_2A @AAd_2A 
\end{CD}\\
\hbox to 11 cm{\dotfill}
\end{gather*}


Let us note that the second row is exactly the Hochschild complex, and the
first column is the bar-complex.

Let us introduce the grading on this bicomplex as follows:
\begin{equation}
\Hom^i(A)=\bigl\{A^{\otimes k}\to A^{\otimes l},\quad k-l=i\bigr\}.
\end{equation}
It is clear that $\Hom^\ndot$ is a dg Lie algebra with the bracket constructed
in Section~1.

\subsection{}

\begin{theorem}
For any associative algebra $A$ with unit the cohomology of the complex
$\Hom^\ndot(A)$ \emph(see formula~\emph{(6))} is equal to $\C[0]$.
\end{theorem}

\begin{proof}
It is clear that the complex $\Hom^\ndot(A)$ is equal to the Hochschild
complex $\Hoch^\ndot(A,B^\ndot)$ where
$$
B^\ndot=0\longleftarrow\C\overset0\longleftarrow
A\overset{d_2}\longleftarrow A^{\otimes2}\overset{d_2}\longleftarrow
A^{\otimes3}\longleftarrow\dots
$$
considered as a complex of $A$-bimodules, and $\C$ is considered as
$A$-bimodule with zero action. On the other hand, $B^\ndot$ is the
bar-complex, it is \emph{quasi-isomorphic} to $\C[0]$ for any algebra~$A$
with unit as complex of bimodules, where $\C[0]$ is equipped with
\emph{zero} action.
Then, $\Hoch^\ndot(A,B^\ndot)\overset{qis}\simeq\Hoch^\ndot(A,\C)$. The last
complex is dual to the bar-complex.
\end{proof}

\subsection{}
\begin{theorem}
Let $A=S^\ndot(V)_0$ be the algebra of polynomials on a vector space~$V$
without unit. Then
$$
H^\ndot(\Hom^\ndot(A))\simeq\nw^\ndot(V)\otimes\nw^\ndot(V^*).
$$
\end{theorem}

\begin{proof}
The proof is analogous to the proof of theorem in Section~3, but now the
bar-complex~$B^\ndot$ quasi-isomorphic to $\nw^\ndot(V)$ as a complex of
$A$-bimodules, when $\nw^\ndot(V)$ is equipped with \emph{zero}
structure of $S^\ndot(V)_0$-bimodule.
\end{proof}

\begin{theorem}
The induced Lie algebra structure on
$\nw^\ndot(V)\otimes\nw^\ndot(V^*)$ is the structure arising from
the Clifford algebra $\Clif(V\oplus V^*)$ \emph(with the bracket $[a,b]=a\cdot
b\pm b\cdot a$\emph).
\end{theorem}

\qed

Then $\nw^\ndot(V)\otimes\nw^\ndot(V^*)$ is a graded Lie algebra
$\Clif^\ndot(V\oplus V^*)$, where
$$
\Clif^i(V\oplus V^*)=\bigl\{\Hom(V^{\otimes k},V^{\otimes l}),\
k-l=i\bigr\}.
$$

\subsection{Formality conjecture for Clifford algebra}

\begin{conjecture}
The differential graded Lie algebra $\Hom^\ndot(S(V)_0)$ is quasi-isomorphic
to the Clifford Lie algebra $\Clif^\ndot(V\oplus V^*)$.
\end{conjecture}

It seems that the Hochschild cohomology $H\!H^\ndot(S(V)_0)$ are equal to
polyvector fields vanishing at zero. 
Then
the Hochschild complex $\Hoch^\ndot(S(V)_0)$ is \emph{not} formal as Lie
algebra (the graph

\sevafig{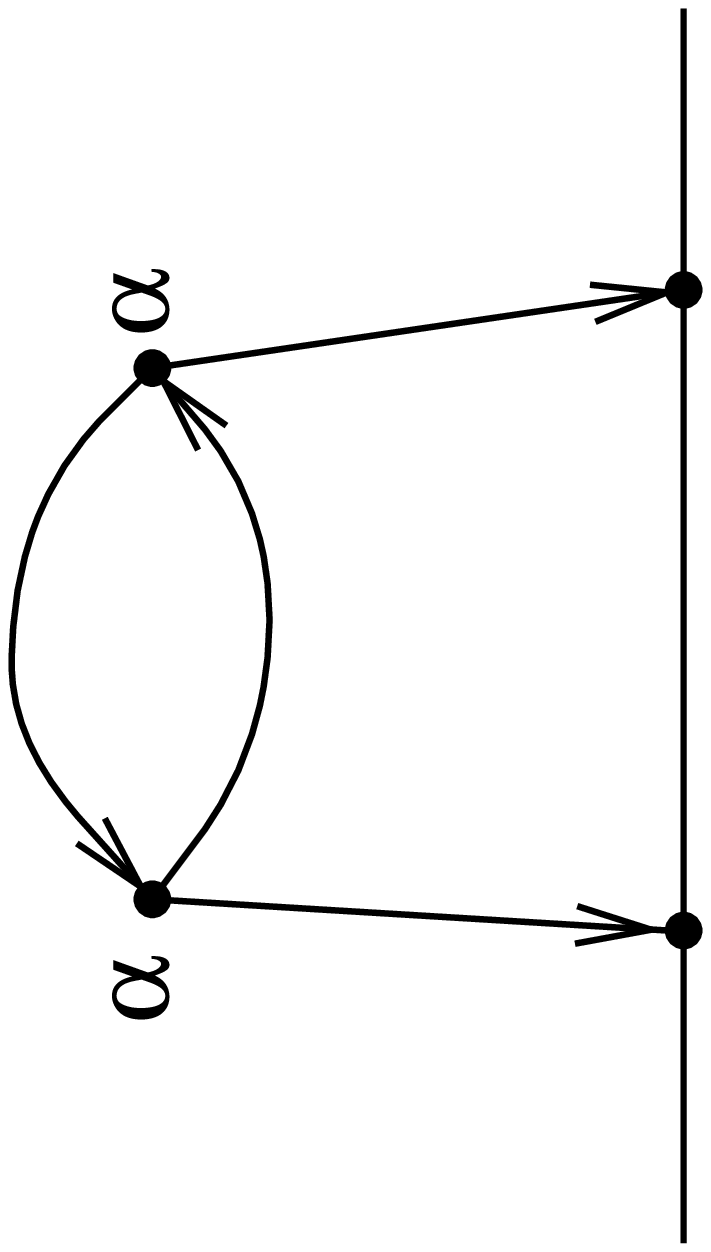}{30mm}{270}

defines a map $S(V)_0^2\to S(V)$, where $\alpha$ is a linear polyvector
field, see~[K]).

As  a  consequence,  there  may  exist   \emph{many}   deformation 
quantizations on
the algebra $S(V)_0$, corresponding to a Lie algebra 
structure on~$V$ or, more generally, to a Poisson bivector field
vanishing at~$0$, which define \emph{the same} (gauge equivalent)
star-products on the algebra~$S(V)$. But any star-product on~$S(V)_0$, 
i.\,e.\ a map $\Psi\colon S(V)_0^{\otimes2}\to
S(V)_0$ (satisfying the Maurer--Cartan equation), 
defines a solution of the Maurer--Cartan equation in the dg Lie algebra
$\Hom^\ndot(S(V)_0)$ (because $d_2=0$ at this place). Then, if the
Conjecture is true, it defines a solution of the Maurer--Cartan equation in
$\Clif(V\oplus V^*)$ modulo the guage equivalence, i.\,e.\ a map
$$
Q\colon \nw^\ndot(V)\to\nw^{\ndot-1}(V)
$$
of degree $-1$ such that $Q^2=0$,
modulo the action
$Q\to Q'=AQA^{-1}$, where $A$ is a map $\nw^\ndot(V)\to\nw^\ndot(V)$
preserving the grading. It is exactly a structure of complex on the graded
space $\nw^\ndot V$
$$
0\longrightarrow\C\overset{Q_0}\longrightarrow
V^*\overset{Q_1}\longrightarrow\nw^2V^*
\overset{Q_2}\longrightarrow\nw^3V^*\longrightarrow\ldots\overset{Q_{n-1}}\longrightarrow\nw^nV^*\longrightarrow0,
$$
i.\,e.\ $Q_iQ_{i+1}=0$,
modulo changes $Q_i\to A_{i+1}Q_iA_i^{-1}$ for linear automorphisms
$(A_0,\dots,A_n)$,\ \ $A_i\in\Aut\nw^iV$. Such a data is described by
dimensions of the cohomology spaces, i.\,e.\ by numbers $b_0,\dots,b_n$ such
that 
\begin{equation}
0\le b_i\le\binom in,\quad \text{and}\quad
\sum_{i=0}^n(-1)^ib_i=\sum_{i=0}^n(-1)^i\binom in=0.
\end{equation}

In the case when $Q=\sum c_{ij}^kv_k\otimes v_i^*\otimes v_j^*\in\nw^2
V^*\otimes V$ is correponded to a Lie algebra structure on~$V$, then
$b_i=\dim H^i_{\Lie}(V;\C)$.

At the moment, I don't know any examples of star-products which would be
interesting from this viewpoint. 

Thus, a star-product on $S^\ndot(V)_0$ is
described (may be not completely) by the corresponding bivector field and by
a sequence $(b_i)$ satisfying~(7).

\subsection{}

In any case, the dg Lie algebra structure on $\Hom^\ndot(A)$ (for any
associative algebra~$A$) produces a quite strange definition of a ``homotopy
algebra structure'' on~$A$ (via the Maurer--Cartan equation). The theorem of
Section~3 shows that for algebras with unit this structure is nondeformable; in
particular, any star-product is equal to usual commutative product in this
sense. But the algebras without unit may give us some interesting examples.

 \end{document}